\documentclass[twoside, 11pt]{article}
\usepackage{mathrsfs,amsfonts,amsmath}
\RequirePackage{ifthen,calc}
\usepackage{theorem}
\usepackage{color}

\usepackage[colorlinks]{hyperref}

\def\R{{\mathbb R}}

\def\E{{\mathbb E}}
\def\P{{\mathbb P}}

\def\N{{\mathbb N}}

\def\qed{\hfill$\square$\smallskip}

 \catcode`@=11
 \def\@evenhead{\hbox to\textwidth{\footnotesize\rm\thepage \hfill
  {\it }}} % authors name

 \def\@oddhead{\hbox to \textwidth{\footnotesize{\it
 } \hfill\thepage}}% abbreviate title

 \renewcommand{\section}{\makeatletter
 \renewcommand{\@seccntformat}[1]{{\csname the##1\endcsname.}\hspace{0.45em}}
 \makeatother \@startsection
{section}%                                            the name
{1}%                                                  the level
{0pt}%                                                the indent
{\baselineskip}%                                      the beforeskip
{0.5\baselineskip}%                                   the afterskip
{\normalsize\bfseries\mathversion{bold}}}

\renewcommand{\subsection}{\makeatletter
 \renewcommand{\@seccntformat}[1]{{\csname the##1\endcsname.}\hspace{0.45em}}
 \makeatother \@startsection
{subsection}%                                            the name
{1}%                                                  the level
{0pt}%                                                the indent
{\baselineskip}%                                      the beforeskip
{0.5\baselineskip}%                                   the afterskip
{\normalsize\bfseries\mathversion{bold}}}
\catcode`@=12

\newtheorem{theorem}{\noindent Theorem}[section]
\newtheorem{lem}{\noindent Lemma}[section]
\newtheorem{cor}{\noindent Corollary}[section]

\theoremheaderfont{\normalfont\bfseries}
\theorembodyfont{\slshape} {\theorembodyfont{\rmfamily}} {\theorembodyfont{\rmfamily}
\newtheorem{defn}{\noindent Definition}[section]}
{\theorembodyfont{\rmfamily} } {\theorembodyfont{\rmfamily}
}
{\theorembodyfont{\rmfamily} } {\theorembodyfont{\rmfamily}
\newtheorem{rem}{\noindent Remark}[section]}
{\theorembodyfont{\rmfamily} } {\theorembodyfont{\rmfamily} } {\theorembodyfont{\rmfamily}
}
\def\beqlb{\begin{eqnarray}}\def\eeqlb{\end{eqnarray}}
\def\beqnn{\begin{eqnarray*}}\def\eeqnn{\end{eqnarray*}}
\numberwithin{equation}{section}
\begin{document}
\title{\bf  Operator Fractional Brownian Motion and Martingale Differences
\footnotetext{\hspace{-5ex}
${[a]}$ College of Mathematics and Information Sciences, Guangxi University,
Nanning, 530004, P.R. China
${[b]}$
Department of Mathematics,
National Tsing-Hua University
Hsinchu, Taiwan, 30043
 \newline
}}
\author{Hongshuai Dai$^{[a]}$,
Tien-Chung  Hu$^{[b]}$, June-Yung Lee$^{[b]}$}

 \maketitle

\begin{abstract}
\noindent  It is well known that martingale difference sequences are very useful in applications and theory. On the other hand, the operator fractional  Brownian motion as an extension of the well-known fractional Brownian motion also plays  important role in both applications and theory. In this paper, we study the relationship between them.  We will construct an approximation sequence of operator fractional  Brownian motion based on a martingale difference sequence.



\end{abstract}

\bigskip

\section{Introduction}
Fractional Brownian motion (FBM) is a continuous Gaussian process with stationary increments. It is one of the well-known self-similar processes. Some studies of financial time series and telecommunication networks have shown that this kind of process with long-range dependency-memory might be a better model in some cases than  the traditional standard Brownian motion. Due to its applications in the real world and its interesting theoretical properties,  fractional Brownian motion has become an object of intense study. One of those studies concerns obtaining its weak limit theorems;  see, for example, Enriquez \cite{E2004}, Niemine \cite{N2004},  Sottinent \cite{S2001}, Li and Dai \cite{LD2011} and the reference therein.

Based on the  study of FBMs, many authors have proposed a generalization of it, and have obtained many new processes. An extension of FBMs is the operator fractional Brownian motion(OFBM). OFBMs are multivariate analogues of one-dimensional FBMs. They arise in the context of multivariate time series and long range dependence (see, for example, Chung \cite{Chung2002}, Davidson
and de Jong \cite{Davidson2000},   Dolado and Marmol \cite{Dolado2004},  Robinson \cite{Robinson2008}, and Marinucci and Robinson \cite{Marinucci2000}). Another context is that of queuing systems, where reflected OFBMs model the size of multiple queues in particular classes of queuing models. They are also studied in problems related to, for example, large deviations (see Delgado \cite{Delgado2007}, and Konstantopoulos and Lin\cite{Lin1996}). Similar to those for FBMs, weak limit theorems for OFBMs have been studied recently. Some new results on approximations of OFBMs have been obtained. See Dai \cite{D2012,D2013} and the references therein.

It is well known that a martingale difference sequence is  extremely useful because it imposes much milder restrictions on the memory of the sequence than under independence, yet most limit theorems that hold for an independent sequence will also hold for a martingale difference sequence.  In recent years,  some researchers have used this type of sequences to construct approximation sequences of some known processes. For example,  Nieminen \cite{N2004} studied the limit theorems for FBMs based on  martingale difference sequences. This is a natural motivation for this paper.  The direct motivation is the recent works by  Dai \cite{D2012, D2013}, in which, based on a sequence of  I.I.D. random variables, the author presented some weak limit theorems for some special kinds of OFBMs.

In this short paper, we establish a weak limit theorem for a special case of OFBMs, which comes from Maejima and Mason~\cite{MM1994}.
The rest of this paper is organized as follows. In Section \ref{Sec2}, we  recall OFBMs and martingale-difference sequences, and present the main result of this paper.  Section \ref{Sec3} is devoted to prove the main result of this paper.

\section{Operator fractional Brownian motion and Martingale-differences}
\label{Sec2}

In this section, we first introduce a special type of OFBMs.
Let $End (\mathbb{R}^{d}) $ be the set of linear operators on $\R^d$
(endomorphisms) and  $Aut (\mathbb{R}^d)$ be the set of invertible linear
operators (automorphisms) in $End (\R^d)$. For convenience, we will
not distinguish an operator $D\in End (\mathbb{R}^d)$ from its associated
matrix relative to the standard basis of $\mathbb{R}^d$. As usual, for $c>0$, $$
c^D=\exp \big((\log c)D\big)=\sum_{k=0}^{\infty} \frac{1}{k!} (\log c)^k D^k.
$$

 Throughout this paper, we will use $\|x\| $ to denote the usual Euclidean norm of $x\in \mathbb{R}^d$. Without confusion, for $A\in End (\mathbb{R}^d)$, we also let
 $\left\|A\right\|=\max_{\|x\|=1}\|Ax\|$ denote the operator norm of $A$. It is easy to see that for $A,B\in End({\mathbb{R}^d})$,
 \beqlb\label{s2-1}
 \left\|AB\right\|\leq \left\|A\right\|\cdot\left\|B\right\|,
 \eeqlb
and for every $A=(A_{ij})_{d\times d}\in End(\mathbb{R}^d)$,
\beqlb\label{s2-2}
\max_{1\leq i,j\leq d}|A_{ij}|\leq \left\|A\right\|\leq d^{\frac{3}{2}} \max_{1\leq i,j\leq d}|A_{ij}|.
\eeqlb
Let $\sigma (A)$ be the collection of all eigenvalues of $A$. We denote
\beqlb\label{s2-3}
\lambda_A=\min\{Re\lambda: \lambda\in\sigma (A)\}\; \hbox{and}\;\Lambda_A=\max\{Re\lambda: \lambda\in\sigma(A)\}.
\eeqlb

 Let $x'$ denote the transpose  of a vector $x\in\mathbb{R}^d$. We now extend the fractional Brownian motion of Riemann-Liouville type studied by L\'{e}vy \cite [p. 357]{levy}  to the multivariate case.

\begin{defn}
Let $D$ be a linear operator on $\mathbb{R}^d$ with
$\frac{1}{2}<\lambda_D,\;\Lambda_D<1$. For $t\in\mathbb{R}_+$, define
\beqlb\label{defn}
X(t)=\int_{0}^t\big(t-u)^{D-I/2} d
W(u),
\eeqlb
where $W(u)=\{W^1(u),\cdots,W^d(u)\}'$ is a standard $d$-dimensional Brownian motion.
We call the process $X=\{X(t)\}$  an operator fractional Brownian motion of Riemann-Liouville (RL-OFBM).
\end{defn}

 As is standard for the multivariate context, we  assume that RL-OFBM is proper.  A random variable in $\mathbb{R}^d$ is  proper if the support of its distribution is not contained in a proper hyperplane of $\mathbb{R}^d$.
\begin{rem}
The operator fractional Brownian motion in the current work is a special case of the operator fractional Brownian  motions in the work of  Maejima and Mason \cite [Theorem 3.1]{MM1994}.
\end{rem}

\begin{rem}\label{es7-lem1}
The RL-OFBM $X$ defined by (\ref{defn}) is an operator self-similar Gaussian process.
\end{rem}

In this short note, we  want to obtain an approximation of RL-OFBMs. Inspired by Nieminen \cite{N2004}, we want to construct an approximation sequence of RL-OFBM $X$ by  martingale differences.

Let $\big\{\xi^{(n)}=(\xi^{(n)}_i,\mathscr{F}_{i}^n)_{1\leq i\leq n}\big\}_{n\in\N}$ be a sequence of square integrable martingale differences such that for every sequence $\{i_n\}$ with $\lim_{n\to\infty}i_n=\infty$, where $1\leq i_n\leq n$,
\beqlb\label{2-1}
\lim_{n\to\infty}\frac{(\xi_{i_n}^{(n)})^2}{\frac{1}{n}}=1,\;\textrm{a.s.}
\eeqlb
and
\beqlb\label{2-2}
\max_{1\leq i\leq n}\big|\xi_i^{(n)}\big|\leq \frac{C}{\sqrt{n}},\;\textrm{a.s.,}
\eeqlb
for some $C\geq 1.$

The following lemma follows from Jacod and Shiryaev \cite{JS1987}.

\begin{lem}\label{lem1}
Under the condition (\ref{2-2}) and the condition
\beqlb\label{2-3}
\sum_{i=1}^{\lfloor nt\rfloor}\big(\xi_i^{(n)}\big)^2\to t,\;\textrm{a.s.,}
\eeqlb
the processes
\beqlb\label{l3-1}
B^n(t)=\sum_{i=1}^{\left\lfloor nt\right\rfloor }\xi_{i}^{(n)}
\eeqlb
converge in distribution to a Brownian motion $B$, as $n\to\infty$.
\end{lem}

\begin{rem}\label{R-1}
Such a type of sequences is very useful, since it is very easy to obtain it in the real world. See, Nieminen \cite{N2004}, for example.
\end{rem}

Below, we extend Lemma \ref{lem1} to the $d$-dimensional case.  Define
\beqlb\label{K-1}
\eta_i^{(n)}=\Big(\xi_{i,1}^{(n)},\cdots,\xi_{i,d}^{(n)}\Big)',
\eeqlb
where $\xi_{i,k}^{(n)}$, $k=1,2\cdots,d$, are independent copies of $\xi_{i}^{(n)}$ in Lemma \ref{lem1}.  Define

\beqlb\label{K-2}
\eta_n(t)=\sum_{i=1}^{\lfloor nt \rfloor}\eta_{i}^{(n)}.
\eeqlb
Then, we can get that $\{\eta^{(n)}\}_{n\in\N}=\{\eta^{(n)}_i,\mathscr{F}^n_i\}$ is still a sequence of square integrable martingale differences on the probability space $\big(\Omega,\mathscr{F}, \P\big)$.  Inspired by Lemma \ref{lem1}, we  have the following lemma.
\begin{lem}\label{lem2}
Under conditions (\ref{2-2}) and (\ref{2-3}), the sequence of processes $\eta_n(t)$ converges in law to a d-dimensional Brownian motion $W$, as $n\to\infty$.
\end{lem}

Noting that $W^i(u)$,~$i=1,\cdots,d$, are mutually independent, and so are  $\xi^{(n)}_{k,i}$,  we can directly get  Lemma \ref{lem2} from Lemma \ref{lem1} and Theorem 11.4.4 in Whitt  \cite[Chapter 12]{W2002}.

Inspired by Lemma \ref{lem2} and (\ref{defn}),  we  construct the approximation sequence by
\beqlb\label{CS}
X_n(t)=\sum_{i=1}^{\left\lfloor nt \right\rfloor}n\int_{\frac{i-1}{n}}^{\frac{i}{n}}\Big(\frac{\left\lfloor nt\right\rfloor}{n}-u\Big)_+^{D-\frac{I}{2}}\eta_{i}^{(n)}du.
\eeqlb

Our main objective in this paper is to explain and prove the following theorem.
\begin{theorem}\label{thm}
The  sequence of processes $\{X_{n}(t),t\in[0,\;1]\}$  given by (\ref{CS}), as $n\to\infty$, converges
weakly to  the operator fractional Brownian motion  $X$  given by (\ref{defn}).
\end{theorem}

In the rest of this paper, most of the estimates  contain unspecified constants. An unspecified positive and finite constant will be denoted by $\tilde{K}$, which may not be the same in each occurrence.

\section {Proof of Theorem \ref{thm}}\label{Sec3}
In order to prove the main result of this paper, we need  a technical lemma.  Before we state this technical lemma, we first introduce the following notation£º
\beqlb\label{defK}
K(t,s)=(t-s)_+^{D-\frac{I}{2}}=\big(K_{i,j}(t,s)\big)_{d\times d},
\eeqlb
and
\beqlb\label{defkn}
K^n(t,s)=\big(\frac{\lfloor nt\rfloor }{n}-s\big)_+^{D-\frac{I}{2}}=\big(K^n_{i,j}(t,s)\big)_{d\times d}.
\eeqlb

The technical lemma follows.
\begin{lem}\label{lem6}
For any $k,j\in\{1,2,\cdots,d\}$,
\beqlb\label{T-1}
\sum_{i=1}^n n^2\int^{\frac{i}{n}}_{\frac{i-1}{n}} K^n_{k,j}(t_l,s)ds\int^{\frac{i}{n}}_{\frac{i-1}{n}} K^n_{k,j}(t_q,s)ds \big(\xi^{(n)}_{i,j}\big)^2\to \int_0^1 K_{k,j}(t_l,s)K_{k,j}(t_q,s)ds,\;\textrm{a.s.}\qquad
\eeqlb
for $t_l,t_q\in[0,\;1]$,  as $n\to\infty$.
\end{lem}

Before we prove it, we need the following lemma which is due to Maejima and Mason \cite{MM1994}.
\begin{lem}\label{lem3}
Let $D\in End(\mathbb{R}^d) $. If $\lambda_D>0$
and $r>0$, then for any $\delta>0$, there exist positive constants $K_1$ and $K_2$ such that
\beqlb\label{s2-4}
\left\|r^D\right\| \leq \begin{cases}K_1 r^{\lambda_D-\delta},  &\textrm{for all}\;  r\leq 1,
\\K_2 r^{\Lambda_D+\delta}, &\textrm{ for all}\; r\geq 1.
\end{cases}
\eeqlb
\end{lem}

Next, we give the detailed proof of Lemma \ref{lem6}.

\noindent{\it Proof of Lemma \ref{lem6}:} In order to simplify the discussion, we split the proof into two steps.

{\bf Step 1.}  We claim that for any $t\in[0,\;1]$,
\beqlb\label{T-2}
\sum_{i=1}^n n^2 \Big(\int_{\frac{i-1}{n}}^{\frac{i}{n}}K_{k,j}(t,s)ds\Big)^2\big(\xi_{i,j}^{(n)}\big)^2\to \int_{0}^1K^2_{k,j}(t,s)ds,\;\textrm{a.s.},
\eeqlb
as $n\to\infty$.

For convenience, define
\beqnn
G_n(t,u)=n \sum_{i=1}^n 1_{(\frac{i-1}{n},\;\frac{i}{n}]}(u)\int_{\frac{i-1}{n}}^{\frac{i}{n}}K_{k,j}(t,s)ds \frac{\xi_{i,j}^{(n)}}{(\sqrt{n})^{-1}}.
\eeqnn
Therefore, we have
\beqlb\label{T-3}
\int_{0}^1 G^2_n(t,u)du&&=\sum_{i=1}^n n^2 \bigg(\int_{\frac{i-1}{n}}^{\frac{i}{n}}K_{k,j}(t,s)ds\bigg)^2 \big(\xi_{i,j}^{(n)}\big)^2\nonumber
\\&&\leq \sum_{i=1}^n n \int_{\frac{i-1}{n}}^{\frac{i}{n}}\big(K_{k,j}(t,s)\big)^2  ds \big(\xi_{i,j}^{(n)}\big)^2,
\eeqlb
where we have used the Cauchy-Schwartz inequality and by \eqref{2-2}.

Therefore,
\beqlb\label{T-39}
\int_{0}^1 G^2_n(t,u)du\leq \tilde{K} \int_0^1 \big( K_{k,j}(t,s)\big)^2ds.
\eeqlb

On the other hand, by (\ref{s2-2}) and Lemma \ref{lem3},
\beqlb\label{T-4}
\big|K_{k,j}(t,s)\big|\leq \big\|K(t,s)\big\|\leq \tilde{K}\big(t-s\big)_+^{(\lambda_D-\delta)-\frac{1}{2}},
\eeqlb
since $t-s\in[0,\;1]$.

By (\ref{T-39}) and (\ref{T-4}), we have
\beqlb\label{T-5}
\int_{0}^1 G^2_n(t,u)du\leq \tilde{K} \int_0^1 \big( K_{k,j}(t,s)\big)^2ds\leq \tilde{K} \int_0^1 \big(t-s\big)_+^{2(\lambda_D-\delta)-1}<\infty,
\eeqlb
since $\lambda_D-\delta>\frac{1}{2}$. Therefore, $\{G^2_n(t,u)\}$ is uniformly integrable.

On the other hand, we have for any $u\in(0,\;1]$,
\beqlb\label{T-6}
G^2_n(t,u)\to K_{k,j}^2(t,u),\;\textrm{a.s.},
\eeqlb
since for $u\in(\frac{i-1}{n},\;\frac{i}{n}]$,
\beqnn
\Big(n\int_{\frac{i-1}{n}}^{\frac{i}{n}}K_{k,j}(t,s)ds\Big)^2\to K^2_{k,j}(t,u),\; \textrm{as}\;n\to\infty,
\eeqnn
and the condition (\ref{2-1}).

By (\ref{T-5}) and (\ref{T-6}), we get that as $n\to\infty$
\beqlb\label{T-7}
\int_{0}^1 G_n^2(t,u)du\to\int_{0}^1 K^2_{k,j}(t,s)ds,\;\textrm{a.s.}
\eeqlb

Therefore, (\ref{T-2}) holds.

{\bf Step 2}. We prove the original claim.  In order to simplify the discussion, we let $t^n_q=\frac{\lfloor nt_q\rfloor}{n}$ and $t^n_l=\frac{\lfloor nt_l\rfloor}{n}$.   By \eqref{T-2},  we can get
\beqlb\label{T-9}
\sum_{i=1}^n n^2\int^{\frac{i}{n}}_{\frac{i-1}{n}} K_{k,j}(t_l,s)ds\int^{\frac{i}{n}}_{\frac{i-1}{n}} K_{k,j}(t_q,s)ds \big(\xi^{(n)}_{i,j}\big)^2\to \int_0^1 K_{k,j}(t_l,s)K_{k,j}(t_q,s)ds,\;\textrm{a.s.}\quad
\eeqlb
for $t_l,t_q\in[0,\;1]$,  as $n\to\infty$.

In fact, it follows from (\ref{T-2}) that
\beqlb\label{T-40}
\sum_{i=1}^n n^2\Big(\int_{\frac{i-1}{n}}^{\frac{i}{n}}K_{k,j}(t_l,s)+K_{k,j}(t_q,s)ds\Big)^2  \big(\xi^{(n)}_{i,j}\big)^2\to \int_0^1\Big(K_{k,j}(t_l,s)+K_{k,j}(t_q,s)\Big)^2 ds.\qquad
\eeqlb
On the other hand, we have
\beqlb\label{T-41}
\Big(\int_{\frac{i-1}{n}}^{\frac{i}{n}}K_{k,j}(t_l,s)+K_{k,j}(t_q,s)ds\Big)^2&&=\Big(\int_{\frac{i-1}{n}}^{\frac{i}{n}}K_{k,j}(t_l,s)ds\Big)^2+\Big(\int_{\frac{i-1}{n}}^{\frac{i}{n}}K_{k,j}(t_q,s)ds\Big)^2\nonumber
\\&&\qquad+2 \int_{\frac{i-1}{n}}^{\frac{i}{n}}K_{k,j}(t_l,s)ds \int_{\frac{i-1}{n}}^{\frac{i}{n}}K_{k,j}(t_q,s)ds.
\eeqlb
Hence  (\ref{T-2}), (\ref{T-40}), and  (\ref{T-41}) imply (\ref{T-9}).

Therefore, in order to prove (\ref{T-1}), it suffices to prove that
\beqlb\label{T-33}
&&\sum_{i=1}^n n^2\bigg(\int^{\frac{i}{n}}_{\frac{i-1}{n}} K_{k,j}(t_l,s)ds\int^{\frac{i}{n}}_{\frac{i-1}{n}} K_{k,j}(t_q,s)ds\nonumber
\\&&\qquad\qquad-\int^{\frac{i}{n}}_{\frac{i-1}{n}} K_{k,j}(t^n_l,s)ds\int^{\frac{i}{n}}_{\frac{i-1}{n}} K_{k,j}(t^n_q,s)ds\bigg)\big(\xi^{(n)}_{i,j}\big)^2 \to 0,\;\textrm{a.s.}
\eeqlb
as $n\to\infty$.

For the left-hand side of (\ref{T-33}), we  have
\beqlb\label{T-34}
&& \int^{\frac{i}{n}}_{\frac{i-1}{n}} K_{k,j}(t_l,s)ds\int^{\frac{i}{n}}_{\frac{i-1}{n}} K_{k,j}(t_q,s)ds -\int^{\frac{i}{n}}_{\frac{i-1}{n}} K_{k,j}(t^n_l,s)ds\int^{\frac{i}{n}}_{\frac{i-1}{n}} K_{k,j}(t^n_q,s)ds\nonumber
\\&&\quad=\int^{\frac{i}{n}}_{\frac{i-1}{n}} K_{k,j}(t_l,s)ds\int^{\frac{i}{n}}_{\frac{i-1}{n}} \Big(K_{k,j}(t_q,s)-K_{k,j}(t^n_q,s)\Big)ds \nonumber
 \\&&\qquad -\int^{\frac{i}{n}}_{\frac{i-1}{n}}\Big( K_{k,j}(t^n_l,s)-K_{k,j}(t_l,s)\Big)ds\int^{\frac{i}{n}}_{\frac{i-1}{n}} \Big(K_{k,j}(t^n_q,s)-K_{k,j}(t_q,s)\Big)ds\nonumber
 \\&&\qquad\quad+\int^{\frac{i}{n}}_{\frac{i-1}{n}} K_{k,j}(t_q,s)ds\int^{\frac{i}{n}}_{\frac{i-1}{n}} \Big(K_{k,j}(t_l,s)-K_{k,j}(t^n_l,s)\Big)ds .
\eeqlb

By (\ref{s2-2}), we have
\beqlb\label{T-35}
 \Big|K_{k,j}(t_q,s)-K_{k,j}(t^n_q,s)\Big|\leq \Big\|K(t_q,s)-K(t_q^n,s)\Big\|.
\eeqlb

On the other hand, using the same method as in the proof of the inequality (\ref{T-22}) below,
\beqlb\label{T-36}
\sum_{i=1}^n \int_{\frac{i-1}{n}}^{\frac{i}{n}} \Big\|K(t_q,s)-K(t_q^n,s)\Big\| ds&&\leq \int^1_0 \Big\|K(t_q,s)-K(t_q^n,s)\Big\|ds\nonumber
\\&&\leq \tilde{K}(t^n_q-t_q)^{2H},
\eeqlb
where $H=\lambda_D-\delta$.

By the condition (\ref{2-2}) and (\ref{T-34}),   (\ref{T-33}) can be bounded by
\beqlb\label{T-37}
&&\tilde{K}n \int^{1}_{0} \Big\|K(t_l,s)\Big\|ds\int^{1}_{0} \Big\|K(t_q,s)-K(t^n_q,s)\Big\|ds \nonumber
 \\&&\qquad + \tilde{K}n \int^{1}_{0}\Big|| K(t^n_l,s)-K(t_l,s)\Big\|ds\int^{1}_{0} \Big\|K(t^n_q,s)-K(t_q,s)\Big\|ds\nonumber
 \\&&\qquad\quad+ \tilde{K}n \int^{1}_{0} \Big\|K(t_q,s)\Big\|ds\int^{1}_{0} \Big\|K(t_l,s)-K(t^n_l,s)\Big\|ds .
\eeqlb

It follows from (\ref{T-5}), (\ref{T-36}), and (\ref{T-37}) that the left-hand side of (\ref{T-33}) can be bounded by
\beqlb\label{T-38}
\tilde{K}n^{1-2H},
\eeqlb
since $|t^n_q-t_q|\leq \frac{1}{n}$ and  $|t^n_l-t_l|\leq \frac{1}{n}$.

From (\ref{T-38}), we can easily prove the lemma. \qed

From the proof of Lemma \ref{lem6} and \eqref{s2-2}, we can easily get that
\begin{cor}\label{lem4}
Let $H(t,s)=\sum_{k=1}^d a_k K^n_{k,j}(t_l,s)$ for any $a_k\in\R$. Then
\beqlb\label{T-56}
\sum_{i=1}^n n^2\int^{\frac{i}{n}}_{\frac{i-1}{n}}H(t_l,s)ds\int^{\frac{i}{n}}_{\frac{i-1}{n}}H(t_q,s)ds \big(\xi^{(n)}_{i,j}\big)^2\to \int_0^1 H(t_l,s)H(t_q,s)ds,\;\textrm{a.s.}
\eeqlb
for any $t_l,t_q\in(0,\;1]$.
\end{cor}
Next, we prove the main result of this paper.  Before we give the details, we first introduce a technical tool.

\begin{lem}\label{a-lem1}
Let $t\in(0,\;1]$, $\sigma_t^2>0$ and let $\{\xi^{(n)}\}$ be a sequence of martingale differences  as in Section \ref{Sec2} and satisfies the following Lindberg condition: for $\epsilon>0$
\beqlb\label{LC}
\sum_{i=1}^{\lfloor nt\rfloor}\E\Big[\big(\xi_i^{(n)}\big)^2 I_{\{|\xi^{(n)}_i|>\epsilon\}}|\mathscr{F}_{i-1}^n\Big]\stackrel{P}{\to}  0.
\eeqlb
Then
\beqlb\label{LC-1}
\sum_{i=1}^{\lfloor nt\rfloor}\big(\xi_i^{(n)}\big)^2 \stackrel{P}{\to} \sigma_t^2
\eeqlb
implies
\beqlb
B^n(t)\stackrel{D}{\to} \mathcal{N}\sim N(0,\;\sigma_t^2),
\eeqlb
where $\stackrel{D}{\to}$ denotes  convergence in distribution.
\end{lem}

Lemma \ref{a-lem1} can be found in  Shiryaev \cite[p. 511]{S1984}.

{\it Proof of Theorem \ref{thm}}:  We will prove this theorem by two steps.

{\bf Step 1}:  First, we have to show that the finite-dimensional distributions of $X_n$ converge to those of $X$.  It suffices to prove that for any $q\in\N$, $a_1,\cdots,a_q\in\R$ and $t_1,\cdots,t_q\in[0,\;1]$,
\beqlb\label{T-8}
\sum_{l=1}^q a_l X_n(t_l)\stackrel{D}{\to}\sum_{l=1}^q a_lX(t_l).
\eeqlb
By the Cram$\acute{e}$r-Wold device (see, Whittle \cite[Chapter 4]{W2002}), in order to prove \eqref{T-8}, we only need to show
\beqlb\label{T-8a}
\sum_{l=1}^q a_l b X_n(t_l)\stackrel{D}{\to}\sum_{l=1}^q a_l b X(t_l),
\eeqlb
for any vector $b=(b^{(1)},\cdots,b^{(d)})\in\R^d$.

For convenience, define
\beqnn
X_n(t)=\Big(X_1^{(n)}(t),\cdots,X_d^{(n)}(t)\Big)',
\eeqnn
where $$
X_j^{(n)}(t)=n\sum_{i=1}^{\lfloor nt \rfloor}\int_{\frac{i-1}{n}}^\frac{i}{n} K_j^n (t,s)\eta^{(n)}_i ds,
$$
with $$K_j^{n}(t,s)=\big(K^n_{j,1}(t,s),\cdots,K^n_{j,d}(t,s)\big),$$
and
\beqnn
X(t)=\Big(X^{(1)}(t),\cdots,X^{(d)}(t)\Big)',
\eeqnn
where
\beqnn
X^{(j)}(t)=\int_{0}^t K_j(t,s)dW(s),
\eeqnn
with
$$K_j(t,s)=\Big(K_{j,1}(t,s),\cdots,K_{j,d}(t,s)\Big).$$
By some calculations, we can get that \eqref{T-8a} is equivalent to
\beqlb\label{T-8aa}
\sum_{l=1}^q\sum_{k=1}^{d}\sum_{j=1}^d\sum_{i=1}^{\lfloor nt_l \rfloor}n\int_{\frac{i-1}{n}}^{\frac{i}{n}}a_l b^{(k)}K^n_{k,j}(t_l,s)\xi_{i,j}^{(n)}ds\stackrel{D}{\to} \sum_{l=1}^q\sum_{k=1}^{d}\sum_{j=1}^{d}\int_{0}^{t}a_l b^{(k)}K_{k,j}(t_l,s)dW^j(s).\qquad
\eeqlb

In order to simplify the discussion, we define
\beqnn
\bar{X}^n(l,k,j)=\sum_{i=1}^{\lfloor nt_l \rfloor}n\int_{\frac{i-1}{n}}^{\frac{i}{n}}K^n_{k,j}(t_l,s)\xi_{i,j}^{(n)}ds
\eeqnn
and
\beqnn
\bar{X}(l,k,j)=\int_{0}^{t_l}K_{k,j}(t_l,s)dW^j(s).
\eeqnn
Hence \eqref{T-8aa} can be rewrote as follows.
\beqlb\label{T-8a3}
\sum_{l=1}^q\sum_{k,j=1}^{d}a_l b^{(k)}\bar{X}^n(l,k,j) \stackrel{D}{\to}\sum_{l=1}^q\sum_{k,j=1}^{d}a_l b^{(k)}\bar{X}(l,k,j).
\eeqlb

By  the independence of  $\xi^{(n)}_{i,j}$, $j=1,\cdots, d$, it suffices to show  that for every  $j\in\{1,\cdots,d\}$
\beqlb\label{T-11a}
\sum_{l=1}^q\sum_{k=1}^d a_lb^{(k)}\bar{X}^n(l,k,j) \stackrel{D}{\to} \sum_{l=1}^q\sum_{k=1}^d a_lb^{(k)} \bar{X}(l,k,j).
\eeqlb
We will prove (\ref{T-11a}) by Lemma \ref{a-lem1}.  We first prove the Lindeberg condition holds in our case. For convenience, define:
\beqnn
Z_{k,i}^n(t)=n\int_{\frac{i-1}{n}}^{\frac{i}{n}}K_{k,j}^n(t,s)\xi_{i,j}^{(n)} ds.
\eeqnn

We have
\beqlb\label{T-13}
\Big(Z_{k,i}^n(t)\Big)^2&&=n^2\big(\xi_{i,j}^{(n)}\big)^2\Big(\int_{\frac{i-1}{n}}^{\frac{i}{n}}K_{k,j}^n(t,s)ds\Big)^2\nonumber
\\&&\leq n \big(\xi_{i,j}^{(n)}\big)^2\int_{\frac{i-1}{n}}^{\frac{i}{n}}\Big(K_{k,j}^n(t,s)\Big)^2ds,
\eeqlb
where we have used the H{\em\"{o}}lder inequality.  By (\ref{s2-2}), we have
\beqlb\label{T-20}
\big|K_{k,j}^n(t,s)\big|\leq \big\|(t-s)_+^{D-\frac{I}{2}}\big\|.
\eeqlb
By Lemma \ref{lem3}, we have
\beqlb\label{T-21}
\big\|(t-s)_+^{D-\frac{I}{2}}\big\|\leq \tilde{K} (t-s)_+^{\lambda_D-\frac{1}{2}-\delta},
\eeqlb
since $t,s\in[0,\;1]$.

By (\ref{T-20}) and (\ref{T-21}),
\beqlb\label{T-43}
\int_{\frac{i-1}{n}}^{\frac{i}{n}}\big\|K_{k,j}(t,s)\big\|^2ds \leq \tilde{K}\int_{\frac{i-1}{n}}^{\frac{i}{n}}(t-s)_+^{2(\lambda_D-\delta)-1}ds\leq \tilde{K}\int_{0}^{\frac{1}{n}}(1-s)^{2(\lambda_D-\delta)-1}ds,
\eeqlb
since $(1-s)^{2(\lambda_D-\delta)-1}$  with $\lambda_D-\delta>\frac{1}{2}$ is decreasing in $s$.  Therefore, it follows from \eqref{T-13} and \eqref{T-43} that
\beqlb\label{T-49}
\Big(Z_{k,i}^n(t)\Big)^2\leq \tilde{K}n \big(\xi_{i,j}^{(n)}\big)^2 \delta_n
\eeqlb
with $\delta_n=\int_{0}^{\frac{1}{n}}(1-s)^{2(\lambda_D-\delta)-1}ds$.

On the other hand, from  \eqref{defkn}, we get for any $s\geq\frac{\lfloor nt \rfloor}{n}$,
\beqlb\label{T-45}
K^n_{k,j}(t,s)=0.
\eeqlb
Hence, by \eqref{T-45},
\beqlb\label{T-46}
\sum_{l=1}^q \sum_{k=1}^d a_l b^{(k)} \bar{X}^n(l,k,j)=\sum_{i=1}^n\sum_{l=1}^q \sum_{k=1}^d a_l b^{(k)}Z_{k,i}^n(t_l).
\eeqlb
Finally, we have
\beqlb\label{T-47}
\big(\sum_{l=1}^q \sum_{k=1}^d a_l b^{(k)}Z_{k,i}^n(t_l)\big)^2&&\leq \tilde{K}\sum_{l=1}^q \sum_{k=1}^d (b^{(k)})^2 a_l^2\big(Z_{k,i}^n(t_l)\big)^2.
\eeqlb
Combining  \eqref{T-43} and \eqref{T-47}, we have
\beqlb\label{T-48}
\big(\sum_{l=1}^q\sum_{k=1}^d a_lb^{(k)}Z_{k,i}^n(t_l)\big)^2\leq \tilde{K}n\big(\xi_{i,j}^n\big)^2\delta_n.
\eeqlb

Noting that
\beqlb\label{T-12}
\Big\{|\sum_{l=1}^q \sum_{k=1}^da_lb^{(k)}Z_{k,i}^n(t_l)|>\epsilon \Big\}=\Big\{\big(\sum_{l=1}^q\sum_{k=1} ^d a_lb^{(k)}Z_{k,i}^n(t_l)\big)^2>\epsilon^2\Big\},
\eeqlb from \eqref{T-48}, we have
\beqlb\label{T-14}
\Big\{|\sum_{l=1}^q \sum_{k=1}^d a_l b^{(k)}Z_{k,i}^n(t_l)|>\epsilon\Big\}\subset \Big\{\tilde{K}n\big(\xi_{i,k}^n\big)^2\delta_n>\epsilon^2 \Big\}.
\eeqlb

Therefore, by  \eqref{T-48} and \eqref{T-14}
\beqlb\label{T-15}
&&\E\Big(\big(\sum_{l=1}^q \sum_{k=1}^d a_lb^{(k)}Z_{k,i}^n(t_l)\big)^2I_{\{|\sum_{l=1}^q \sum_{k=1}^da_lb^{(k)}Z_{k,i}^n(t_l)|>\xi\}}|\mathscr{F}_{i-1}^n\Big)\nonumber
\\&&\qquad\qquad\leq \tilde{K}n\big(\xi_{i,j}^n\big)^2\delta_n \E\Big(I_{\Big\{\tilde{K}n\big(\xi_{i,j}^n\big)^2\delta_n>\epsilon^2 \Big\}}|\mathscr{F}_{i-1}^n\Big)\nonumber
\\ &&\qquad\qquad \leq \tilde{K} \delta_n \E\Big(I_{\Big\{\tilde{K}\delta_n>\epsilon^2 \Big\}}|\mathscr{F}_{i-1}^n\Big).
\eeqlb

Combining \eqref{T-46} and \eqref{T-15}, one can easily prove that, as $n$ approaches $\infty$,
\beqnn
\sum_{i=1}^n \E\Big(\big(\sum_{l=1}^{q}\sum_{k=1}^d a_lb^{(k)}Z_{k,i}^n\big)^2I_{\{|\sum_{l=1}^{q}\sum_{k=1}^da_lb^{(k)}Z_{k,i}^n|>\xi\}}|\mathscr{F}_{i-1}^n\Big)\to 0.
\eeqnn
Hence the Lindeberg condition holds.

Next, we show the condition (\ref{LC-1}) holds. We first study the right-hand side of \eqref{T-11a}.  We have
\beqlb\label{T-50}
\sum_{l=1}^q\sum_{k=1}^{d}a_l b^{(k)}\bar{X}(l,k,j)=\sum_{l=1}^q a_l \tilde{W}(t_l),
\eeqlb
where \beqlb\label{T-51}
\tilde{W}(t)=\int_0^t \Big[\sum_{k=1}^d b^{(k)}K_{k,j}(t,s)\Big]dW^j(s)= \int_0^t \bar{K}(t,s)dW^j(s),
\eeqlb
with
$$
\bar{K}(t,s)=\sum_{k=1}^d b^{(k)}K_{k,j}(t,s).
$$
Combining \eqref{T-50} and \eqref{T-51}, we have
\beqlb\label{T-52}
\E\Big[\sum_{l=1}^q\sum_{k=1}^{d}a_l b^{(k)}\bar{X}(l,k,j)\Big]^2=\E\Big[\sum_{l=1}^q a_l \tilde{W}(t_l)\Big]^2=\sum_{l,j=1}^q a_ja_l\int_0^1\bar{K}(t_j,s)\bar{K}(t_l,s)ds.\qquad
\eeqlb
In order to show the condition (\ref{LC-1}), we only need to show
\beqlb\label{T-53}
\sum_{i=1}^n \big(\sum_{l=1}^q \sum_{k=1}^d a_l b^{(k)}Z_{k,i}^n\big)^2\stackrel{P}{\to} \sum_{l,j=1}^q a_ja_l\int_0^1\bar{K}(t_j,s)\bar{K}(t_l,s)ds.
\eeqlb

Now, we focus on the left-hand side of \eqref{T-53}.  Similar to \eqref{T-50}, we have
\beqlb\label{T-54}
\sum_{l=1}^q \sum_{k=1}^d a_l b^{(k)}Z_{k,i}^n=\sum_{l=1}^q a_l \bar{Z}_{l,i}^n,
\eeqlb
where
\beqlb\label{T-55}
\bar{Z}_{l,i}^n=n\int_{\frac{i-1}{n}}^{\frac{i}{n}}\bar{K}^n_j(t_l,s)\xi^{(n)}_{i,j}ds
\eeqlb
with $\bar{K}^n_j(t_l,s)=\sum_{k=1}^d b^{(k)}K_{k,j}^n(t_l,s)$.  Hence
\beqlb\label{T-16}
\sum_{i=1}^n \big(\sum_{l=1}^q \sum_{k=1}^d a_l b^{(k)}Z_{k,i}^n\big)^2&&=\sum_{i=1}^n \sum_{l_1,l_2=1}^q n^2 a_{l_1}a_{l_2} \int_{\frac{i-1}{n}}^{\frac{i}{n}}\bar{K}_{j}(t_{l_1},s)ds\nonumber
\\&&\qquad\int_{\frac{i-1}{n}}^{\frac{i}{n}}\bar{K}_{j}(t_{l_2},s)ds\big(\xi^{(n)}_{i,j}\big)^2.
\eeqlb
It follows from Corollary \ref{lem4} that the right-hand side of the equation \eqref{T-16} converges to
\beqlb\label{T-17}
\sum_{l_1,l_2=1}^qa_{l_1}a_{l_2}\int_0^1 \bar{K}_{j}(t_{l_1},s)\bar{K}_{j}(t_{l_2},s)  ds, \;\textrm{a.s.}
\eeqlb
as $n\to\infty$.
On the other hand, one can easily get that
\beqlb\label{T-17a}
\E\Big[\tilde{W}(t_l)\tilde{W}(t_k)\Big]=\int_0^1 \bar{K}(t_l,s)\bar{K}(t_k,s)ds.
\eeqlb

By \eqref{T-52}, \eqref{T-17}, and (\ref{T-17a}), we get the condition (\ref{LC-1}).

{\bf Step 2}: We need to prove the tightness of the sequence $\{X_n(t)\}$.

By some calculations,
\beqlb\label{T-19}
\E\Big(\|X_n(t)-X_n(s)\|^2\Big)\leq \tilde{K} \int_{0}^1 \Big\|\Big(\frac{\lfloor nt\rfloor}{n}-u\Big)_+^{D-\frac{I}{2}}-\Big(\frac{\lfloor ns\rfloor}{n}-u\Big)_+^{D-\frac{I}{2}}\Big\|^2 du.
\eeqlb
In order to simplify the discussion, let
\beqnn
\tilde{t}=\frac{\lfloor nt\rfloor}{n},\;\textrm{and}\;\tilde{s}=\frac{\lfloor ns\rfloor}{n}.
\eeqnn
Next, we show that
\beqlb\label{T-22}
\int_{0}^1 \Big\|\Big(\tilde{t}-u\Big)_+^{D-\frac{I}{2}}-\Big(\tilde{s}-u\Big)_+^{D-\frac{I}{2}}\Big\|^2 du\leq K(\tilde{t}-\tilde{s})^{2H},
\eeqlb
where $H=\lambda_D-\delta$.

In fact,
\beqlb\label{T-23}
&&\int_0^1 \Big\|(\tilde{t}-u)_+^{D-\frac{I}{2}}-(\tilde{s}-u)_+^{D-\frac{I}{2}}\Big\|^2du\nonumber
\\&&\qquad\qquad =\int_0^{\tilde{s}} \big\|(\tilde{t}-u)^{D-\frac{I}{2}}-(\tilde{s}-u)^{D-\frac{I}{2}}\big\|^2 du\nonumber
\\&&\qquad\qquad\quad+\int_{\tilde{s}}^{\tilde{t}} \big\|(\tilde{t}-u)^{D-\frac{I}{2}}\big\|^2du.
\eeqlb

It follows from Lemma \ref{lem3} and  (\ref{s2-1}) that
\beqnn
\big\|(\tilde{t}-u)^{D-\frac{I}{2}}\big\| \leq \tilde{K} (\tilde{t}-u)^{\lambda_D-\delta-\frac{1}{2}},
\eeqnn
since $ u\leq \tilde{t}\in[0,1]$.

Therefore,
\beqlb\label{T-24}
\int_{\tilde{s}}^{\tilde{t}} \big\|(\tilde{t}-u)_+^{D-\frac{I}{2}}\big\|^2du &&\leq \tilde{K} \int_{\tilde{s}}^{\tilde{t}}(\tilde{t}-u)^{2(\lambda_D-\delta)-1}du\nonumber
\\&&=\frac{\tilde{K}(\tilde{t}-\tilde{s})^{2(\lambda_D-\delta)}}{2(\lambda_D-\delta)}.
\eeqlb

Next, we deal with the first term on the right-hand side of (\ref{T-23}).  Note that
\beqlb\label{T-25}
&&\int_0^{\tilde{s}} \big\|(\tilde{t}-u)^{D-\frac{I}{2}}-(\tilde{s}-u)^{D-\frac{I}{2}}\big\|^2 du=\int_0^{\tilde{s}} \big\|(\tilde{t}-\tilde{s}+u)^{D-\frac{I}{2}}-u^{D-\frac{I}{2}}\big\|^2 du\nonumber
\\&&\qquad\qquad= \int_0^{\tilde{s}/(\tilde{t}-\tilde{s})} \Big\|\big[(\tilde{t}-\tilde{s})(1+u)\big]^{D-\frac{I}{2}}-\big[(\tilde{t}-\tilde{s})u\big]^{D-\frac{I}{2}}\Big\|^2 du (\tilde{t}-\tilde{s})\nonumber
\\&&\qquad\qquad\leq \big\|(\tilde{t}-\tilde{s})^{D-\frac{I}{2}}\big\|^2 (\tilde{t}-\tilde{s}) \int_0^{\tilde{s}/(\tilde{t}-\tilde{s})} \big\|(1+u)^{D-\frac{I}{2}}-u^{D-\frac{I}{2}}\big\|^2 du\nonumber
\\&&\qquad\qquad \leq \big\|(\tilde{t}-\tilde{s})^{D-\frac{I}{2}}\big\|^2 (\tilde{t}-\tilde{s})\int_{\mathbb{R}_+} \Big\|(1+u)^{D-\frac{I}{2}}-u^{D-\frac{I}{2}}\Big\|^2 du,
\eeqlb
where we used the fact that $(\tilde{t}\tilde{s})^{A}=\tilde{t}^A\cdot \tilde{s}^A$.

It follows from Lemma \ref{lem3} and  (\ref{s2-1}) that
\beqnn
\big\|(\tilde{t}-\tilde{s})^{D-\frac{I}{2}}\big\|^2 (\tilde{t}-\tilde{s})\leq \tilde{K}(\tilde{t}-\tilde{s})^{2(\lambda_D-\delta)}.
\eeqnn

In order to prove our result, it suffices to show that \beqlb\label{T-26} \int_{\mathbb{R}_+}
\big\|(1+u)^{D-\frac{I}{2}}-u^{D-\frac{I}{2}}\big\|^2 du<\infty. \eeqlb

Then, in order to prove (\ref{T-26}), it suffices to show that
\beqlb\label{T-27} \int_{u\leq 1} \big\|u^{D-\frac{I}{2}}\big\|^2 du
<\infty, \eeqlb and for large enough $T>1$, that \beqlb\label{T-28}
\int_{ u \geq T} \big\|(1+u)^{D-\frac{I}{2}}-u^{D-\frac{I}{2}}\big\|^2
du<\infty. \eeqlb

It follows from Lemma \ref{lem3} and  (\ref{s2-1}) that \beqnn
\big\|u^{D-\frac{I}{2}}\big\|^2\leq \tilde{K}u^{2(\lambda_D-\delta)-1} \;\textrm
{for}\; u\leq 1. \eeqnn Hence, one can easily see that
(\ref{T-26}) holds.

Next, we show that (\ref{T-28}) holds. We  see that \beqlb\label{es8-15}
(1+u)^{D-\frac{I}{2}}-u^{D-\frac{I}{2}}=\int_u^{1+u}(D-\frac{I}{2})s^{D-\frac{I}{2}}
s^{-1}ds. \eeqlb Then \beqlb\label{es8-16}
\big\|(1+u)^{D-\frac{I}{2}}-u^{D-\frac{I}{2}}\big\| \leq
\big\|(D-\frac{I}{2})\big\| \int_u^{1+u}||s^{D-\frac{I}{2}}\big\| s^{-1}ds.
\eeqlb It follows from Lemma \ref{lem1} and  (\ref{s2-1}) that
\beqlb\label{es8-17}
 \int_u^{1+u}\big\|s^{D-\frac{I}{2}}\big\| s^{-1}ds\leq \int_u^{1+u} \tilde{K} s^{\Lambda_D+\delta-\frac{3}{2}}ds,
\eeqlb
since $u\geq 1$.

By (\ref{es8-16}) and (\ref{es8-17}),
\beqlb\label{es8-18}
\big\|(1+u)^{D-\frac{I}{2}}-(u)^{D-\frac{I}{2}}\big\|^2 \leq \tilde{K} u^{2(\Lambda_D+\delta)-3}.
\eeqlb

 By (\ref{es8-18}), we  have that (\ref{T-28}) holds, since $\Lambda_D+\delta<1$

Therefore, we have
\beqlb\label{T-29}
\E\Big(\|X_n(t)-X_n(s)\|^2\Big)\leq \tilde{K}(\tilde{t}-\tilde{s})^{2H}.
\eeqlb
 Hence for any $s\leq t\leq u\in [0,\;1]$, we have
\beqlb\label{T-42}
\E\bigg[\Big\|X_n(t)-X_n(s)\Big\|\Big\|X_n(t)-X_n(u)\Big\|\bigg]&&\leq \Big[\E\big\|X_n(t)-X_n(s)\big\|^2\Big]^{\frac{1}{2}}\Big[\E\big\|X_n(t)-X_n(u)\big\|^2\Big]^{\frac{1}{2}}\nonumber
\\ &&\leq \tilde{K}\Big|\frac{\lfloor nt \rfloor}{n}-\frac{\lfloor ns \rfloor}{n} \Big|^{H}\Big|\frac{\lfloor nu \rfloor}{n}-\frac{\lfloor nt \rfloor}{n} \Big|^{H}\nonumber
\\ &&\leq \tilde{K} \Big|\frac{\lfloor nu \rfloor}{n}-\frac{\lfloor ns \rfloor}{n} \Big|^{2H}.
\eeqlb

If $u-s\geq \frac{1}{n}$, then one can easily see that
\beqlb\label{T-31}
\E\bigg[\Big\|X_n(t)-X_n(s)\Big\|\Big\|X_n(t)-X_n(u)\Big\|\bigg]\leq \tilde{K} (u-s)^{2H}.
\eeqlb
On the other hand, if  $u-s<\frac{1}{n}$, then either $s$ and $t$ or $t$ and $u$ belong to the interval $[\frac{i}{n},\;\frac{i+1}{n}]$ for some $i$. Thus the left-hand side of (\ref{T-42}) is zero.  Therefore (\ref{T-31}) still holds for this case.  Hence it follows  from Ethier and Kurtz \cite[Chapter 3]{EK86} that $\{X_n(t)\}$ is tight.

By Theorem 7.8 in Ethier and Kurtz \cite[Chapter 3]{EK86}, we  get that Theorem \ref{thm} holds.  This completes the proof.\qed

\medskip
\noindent{\bf Acknowledgments}\ This work was  supported by the  Natural Science Foundation of  China (No. 11361007) and the Guangxi Natural Science Foundation (No. 2012NSFGXBA05301).

\end{document}